\documentclass[11pt]{article}

\usepackage{amsmath}
\usepackage{amsmath,amsfonts,amssymb,amsthm,graphics,amscd}
\usepackage[lmargin=3cm, rmargin=3cm,tmargin=2.5cm,bmargin=2.5cm]{geometry}

\usepackage{color}

\newcommand{\bdf}{\begin{definition}}
\newcommand{\edf}{\end{definition}}
\newcommand{\bprop}{\begin{proposition}}
\newcommand{\eprop}{\end{proposition}}
\newcommand{\bcor}{\begin{corollary}}
\newcommand{\ecor}{\end{corollary}}
\newcommand{\bet}{\begin{theorem}}
\newcommand{\eet}{\end{theorem}}
\newcommand{\blm}{\begin{lemma}}
\newcommand{\elm}{\end{lemma}}
\newcommand{\bp}{\begin{proof}}
\newcommand{\ep}{\end{proof}}
\newcommand{\bex}{\begin{example}\rm}
\newcommand{\eex}{\end{example}}
\newcommand{\bexs}{\begin{examples}\rm}
\newcommand{\eexs}{\end{examples}}
\newcommand{\bremark}{\begin{remark}\rm}
\newcommand{\eremark}{\end{remark}}

\newcommand{\bnotations}{\begin{notations}\rm}
\newcommand{\enotations}{\end{notations}}

\newcommand{\norm}[1]{ \| #1 \| }

\newtheorem{theorem}{Theorem}[section]
\newtheorem{lemma}[theorem]{Lemma}
\newtheorem{corollary}[theorem]{Corollary}

\newtheorem{example}[theorem]{Example}
\newtheorem{definition}[theorem]{Definition}
\newtheorem{proposition}[theorem]{Proposition}
\newtheorem{remark}[theorem]{Remark}

\newtheorem{examples}[theorem]{Examples}
\newtheorem{notations}[theorem]{Notations}

\numberwithin{equation}{section}

\begin{document}

\title {Fredholm-type  Operators and Index }

\author { Alaa Hamdan and  Mohammed
Berkani\ }

\date{}

\maketitle \setcounter{page}{1}

\begin{abstract} While in \cite{HB}  we studied classes of Fredholm-type operators defined by the homomorphism $\Pi$ from  $L(X)$ onto the Calkin algebra  $\mathcal{C}(X)$, $X$ being  a Banach space, we study  in this paper two classes of Fredholm-type operators defined by the homomorphism $\pi$ from  $L(X)$ onto the algebra  $\mathcal{C}_0(X)= L(X)/F_0(X),$  where $F_0(X)$ is the ideal of finite rank operators in $L(X).$   Then we   define  an index   for  Fredholm-type operators and we show  that this new index satisfies  similar properties as the usual Fredholm index.

 \end{abstract}

\renewcommand{\thefootnote}{}
\footnotetext{\hspace{-7pt} {\em 2020 Mathematics Subject
Classification\/}: 47A53,  16U90, 16U40.
\baselineskip=18pt\newline\indent {\em Key words and phrases\/}:
Calkin algebra,  Fredholm operators,   generalized Drazin invertibility, index, idempotent,
p-invertibility.}

\section{Introduction}

Let $X$ be a Banach space,
let $L(X)$  be the Banach algebra of bounded linear operators acting on
the Banach space $X$ and let  $ T \in L(X).$
We will denote by $N(T)$  the null space of  $T,$ by $ \alpha(T)$
the nullity of $T,$  by $ R(T)$ the range of $T$ and by $\beta(T) $
its defect. If the range $ R(T)$ of $T$ is closed  and  $ \alpha(T) < \infty $
 (resp.  $\beta(T)< \infty $ ),
\noindent  then $T$ is called an upper semi-Fredholm (resp. a
lower semi-Fredholm) operator. A semi-Fredholm operator is an
upper or a lower semi-Fredholm operator. If both of $ \alpha(T) $
and $\beta(T) $ are finite, then $T $ is called a Fredholm operator
and the index of $T$  is defined by  the formula \,  $ ind(T) = \alpha(T) -
\beta(T). $  The notations $\Phi_+(X),  \Phi_-(X)$ and  $\Phi(X)$ will designate  respectively
the set of upper semi-Fredholm, lower semi-Fredholm and Fredholm operators.

%\noindent For $T\in L(X),$ the Fredholm  spectrum of $T$
% is defined by: $$ \sigma_F(T) = \{\lambda \in \CC: T-\lambda I \text { is\ not\ a\ Fredholm\ operator}  \}.$$

\noindent Define also  the sets $\Phi_l(X) = \{T \in \Phi_+ (X) \mid $ there exists a bounded projection of $X$ onto $R(T)$\}  and $\Phi_r (X) = \{T \in \Phi_-(X)\mid $ there exists a bounded projection of $X $ onto $N(T) \}$.

\noindent Recall that the Calkin algebra over $X$ is the Banach algebra given by  the quotient algebra
$\mathcal{C}(X)=L(X)/K(X)$, where  $K(X)$ is the closed ideal of
compact operators on $X$.  Let $G_r$ and $G_l$   be the right and left, respectively, invertible elements
  of $ \mathcal{C}(X)$. From \cite[Theorem 4.3.2]{CPY}    and  \cite[Theorem 4.3.3]{CPY},
   it follows that $ \Phi_l (X)= \Pi^{-1}(G_l)$  and  $ \Phi_r (X)= \Pi^{-1}(G_r),$
     where the homomorphism $ \Pi: L(X) \rightarrow  \mathcal{C}(X)$ is the natural projection. We observe that $\Phi(X)= \Phi_l (X)\cap \Phi_r(X).$

  \bdf The elements of  $ \Phi_l (X)$ and   $\Phi_r (X)$ will be  called respectively left
    semi-Fredholm operators and   right  semi-Fredholm operators.

  \edf

In 1958, in his paper \cite{DR},  M.P. Drazin  extended the
concept of invertibility in  rings and semigroups and introduced a
new kind of inverse, known now as the  Drazin inverse.

\bdf An element  $a$  of  a semigroup $ \mathcal{S}$ is called
Drazin invertible, if there exists an element $b \in \mathcal{S}$  written $ b= a^d$  and called the Drazin inverse of $a,$
satisfying the following equations:
\begin{equation}\label{Drazin}
   ab = ba, b=ab^2, a^k= a^{k+1}b,
 \end{equation}
for some nonnegative integer $k.$ The least nonnegative
integer $k$ for which these equations holds is the Drazin index $
i(a)$ of $a.$

\edf

\noindent It follows from \cite{DR} that a Drazin invertible
element in a semigroup has a unique Drazin inverse.\\
If $A$ is a unital Banach algebra,  it is easily seen that the conditions  (\ref{Drazin}) for the Drazin invertibility of an element $ a\in A$  are equivalent to the existence of an idempotent $ p \in A $ such that $ap=pa, ap \text{ is nilpotent and} \,\,  a+p   \text{ is   invertible}. $

\noindent In 1996, in  \cite[Definition 2.3]{Koliha}, J.J. Koliha extended
the notion of  Drazin invertibility.

\bdf An element  $a$  of a unital Banach algebra $A$ will be said to be
generalized Drazin invertible if there exists $b \in A $ such that
$ bab=b, ab=ba$ and  $ aba-a$ is a quasinilpotent element in $A$.
\edf

 In \cite{Koliha}, Koliha proved that an element  $a\in A$ is generalized Drazin
invertible if and only if  there exists $ \epsilon >0,$ such that
for all complex scalar  $\lambda $ such that $ 0 < \mid \lambda \mid <  \epsilon,
$ the element  $a-\lambda e $ is invertible, where $e$ is the unit of $A.$  He proved in \cite[Theorem 4.2]{Koliha},
 that a generalized Drazin
invertible element  has a unique generalized Drazin
inverse. He also proved that an element  $a\in A$ is generalized Drazin
invertible if and only if  there exists an idempotent $p \in A$  commuting with $a,$
such that $ a+p$ is invertible in $A$ and $ap$ is quasinilpotent. \\
Let $A$ be a ring with a unit and let $p$ be an idempotent element in $A$.  Recall that the  commutant $C_p$ of $p$  is the subring of $A$ defined  by $C_p=
\{ x \in A \mid xp=px \}.$
   In  \cite[Definition 2.2]{P50},
the concepts of left $p-$invertibility, right $p-$invertibility and $p-$invertibility were defined  as follows.

\bdf Let $ a \in A.$   We will say that

\begin{enumerate}
\item $a$ is left  $p$-invertible if $ ap=pa$ and  $ a+p$ is left
invertible in $C_p,$
\item  $a$ is right  $p$-invertible if $ ap=pa$ and $ a+p$ is right
invertible in $C_p,$
 \item $a$ is
$p$-invertible if $ ap=pa$    and $ a+p$ is invertible.
\end{enumerate}

\edf

Moreover in \cite[Definition 2.11]{P50}, left and right Drazin invertibility were defined as follows.

\bdf \label{Drazin $p$-invertibility}  An element
$a  \in A $ is left Drazin invertible (respectively right Drazin
invertible) if there exists an idempotent $p \in A$ such that $a$
is left $p$-invertible (respectively right $p$-invertible) and
$ap$ is nilpotent.

\edf

For  $T \in L(X),$  we will say that  a subspace $M$
 of $X$  is {\em $T$-invariant} if $T(M) \subset M.$  If $M$
  is {\em $T$-invariant,} we define
$T_{\mid M}:M \to M$ as $T_{\mid M}(x)=T(x), \, x \in M$.  If $M$
and $N$ are two closed $T$-invariant subspaces of $X$ such that
$X=M \oplus N$, we say that $T$ is {\em completely reduced} by the
pair $(M,N)$ and it is denoted by $(M,N) \in Red(T)$. In this case,
we write $T=T_{\mid M} \oplus T_{\mid N}$ and we say that $T$ is a
{\em direct sum} of $T_{\mid M}$ and $T_{\mid N}$.
 To say that a pair $(M,N)$ of closed subspaces of $X$ is  in $Red(T),$
means simply that there exists an idempotent $P \in L(X)$ commuting with $T.$

It is said that $T \in L(X)$ admits a  generalized
Kato decomposition, abbreviated as GKD,  if there exists $(M, N)
\in Red(T)$ such that $ T_{\mid M}$ is Kato and $T_{\mid N}$ is
quasinilpotent. Recall that an operator $T \in L(X)$ is {\em Kato}
if $R(T)$ is closed and $Ker(T) \subset R(T^n)$ for every $ n \in
\mathbb{N}$.

\bdf \label{relative-index}

Let $T \in L(X)$ and let $M$ be a closed subspace of $X$ which is  {\em $T$-invariant}.
We define the  nullity of\,  $T$   relatively  to $M$ by $ \alpha_M(T)= \alpha(T_{\mid M}),$  the deficiency of $T$   relatively  to $M$ by $ \beta_M(T)= \beta(T_{\mid M}).$
If $T_{\mid M}$ is a Fredholm operator, then the index  of $T$ relatively to $M$ is defined by the formula \, $ind_M(T)= \alpha_M(T) - \beta_M(T).$

\edf

For $T \in L(X)$ and a nonnegative integer $n,$ define $T_{[n]}$ to
be the restriction of $T$ to $R(T^n)$ viewed as a map from $R(T^n)$
into $R(T^n)$ (in particular $T_{[0]}=T).$ If for some integer $n$
the range space $R(T^n)$ is closed and $T_{[n]}$ is an upper (resp.
a lower) semi-Fredholm operator, then $T$ is called an  upper (resp.
a lower) semi-B-Fredholm operator.
Moreover, if $T_{[n]}$ is a Fredholm operator, then $T$ is called a
B-Fredholm operator (see \cite{P7} for more details). From    \cite[Theorem 2.7 ]{P7}, we know  that   $T \in L(X)$  is  a B-Fredholm operator if
 there exists $(M, N) \in Red(T)$ such that $T_{\mid M}$ is a nilpotent   operator and $T_{\mid N}$ is a Fredholm operator.  The notations $ \Phi_B(X), \Phi_{B^+}(X)$ and  $\Phi_{B^-}(X)$ designates  respectively the set of upper semi-B-Fredholm, lower semi-B-Fredholm and B-Fredholm operators. Then  we have $ \Phi_B(X)= \Phi_{B^+}(X) \cap \Phi_{B^-}(X). $

  From \cite[Theorem 3.4]{P10}, we know that  $T \in L(X)$ is a B-Fredholm operator if and only if $\pi(T)$ is Drazin invertible in the algebra $ L(X)/F_0(X), $ where $F_0(X)$ is the ideal of finite rank operators in $L(X)$ and $\pi: L(X) \rightarrow L(X)/F_0(X)$ is the natural projection. To simplify the notations, we will set  $\mathcal{C}_0(X)= L(X)/F_0(X).$

  Let us recall the following important result, which will be  used throughout this paper.  If $p$ is an idempotent element of $\mathcal{C}(X),$ then from \cite[Lemma 1]{BAR}, we know that there exists an idempotent $P \in L(X)$ such that $ \Pi(P)= p.$

 As this paper is a continuation of \cite{HB}, we begin by giving a brief  summary of the results obtained there. So letting $p=\Pi(P)$ being  an idempotent  in the Calkin algebra and    using the results on $p$-invertibility obtained  in  \cite{P50},   we introduced  the class $\Phi_P(X)$ of $P$-Fredholm (respectively the class $\Phi_{lP}(X)$ of left semi-$P$-Fredholm and the class $\Phi_{rP}(X)$ of right  semi-$P$-Fredholm)  operators (Definition \ref{P-Fredholm}), in a similar way as  the corresponding classes of Fredholm, left semi-Fredholm and right semi-Fredholm operators.  We noticed that
$\Phi_P(X)= \Phi_{lP}(X)\cap \Phi_{rP}(X).$

 Based on left and right  Drazin invertibility in the Calkin algebra (see Definition \ref{Drazin $p$-invertibility}), we
 introduced  the classes $\Phi_{l\mathcal{W}B}(X),$  $\Phi_{r\mathcal{W}B}(X) $  of  left and right  weak semi-B-Fredholm operators(\cite [Definition 4.3]{HB}), completing in this way the study of the class  $\Phi_{\mathcal{W}B}(X) $ of weak B-Fredholm operators begun   in \cite{P45} and we proved that  $$\Phi_{\mathcal{W}B}(X)= \Phi_{l\mathcal{W}B}(X)\cap \Phi_{r\mathcal{W}B}(X).$$

 Similarly, using left and right generalized Drazin invertibility in the Calkin algebra (See Definition \ref{gdrazin}), we  introduced the classes $\Phi_{l\mathcal{P}B}(X) $ and $\Phi_{r\mathcal{P}B}(X) $ of left and right  pseudo semi-B-Fredholm operators (See Definition \ref{pseudo}), completing in this way the study of the class  $\Phi_{\mathcal{PB}}(X) $ of pseudo B-Fredholm operators inaugurated  in \cite{P45} and we proved that   \mbox{$\Phi_{\mathcal{P}B}(X)= \Phi_{l\mathcal{P}B}(X)\cap \Phi_{r\mathcal{P}B}(X).$}

 In summary, our study in \cite{HB} was based on various type of invertibility in the Calkin algebra  $\mathcal{C}(X)$ and the algebra  homomorphism  $\Pi:L(X)\rightarrow  \mathcal{C}(X).$  So  we  have  defined classes of  Fredholm type operators through the Calkin algebra.  As seen in \cite{HB},  these classes of operators obeys   the following strict inclusions relations

  $$ \Phi(X) \subsetneq \Phi_B(X) \subsetneq \Phi_{\mathcal{W}B}(X) \subsetneq \Phi_{\mathcal{P}B}(X)\subsetneq \Phi_\mathbb{P}(X),$$
   where  $ \Phi_{\mathbb{P}}(X)=\bigcup\limits_{ \{ P\in L(X) \mid P^2=P\}}\Phi_P(X). $

Now, in a similar way to the work done in \cite{HB}, we will study in the second and the third  sections of the  present paper, the classes of  $P$-Fredholm and pseudo-B-Fredholm operators  defined by considering various type of invertibility  in  $\mathcal{C}_0(X)$ and  the  algebra homomorphism $\pi:L(X) \rightarrow   \mathcal{C}_0(X).$ These two  classes of  Fredholm type operators  are so defined through the algebra $\mathcal{C}_0(X).$

 In the fourth section of the present paper, we will associate  an index to left (resp. right) pseudo-semi-B-Fredholm  and pseudo-B-Fredholm operators. This new index satisfies very similar  properties as the  usual  index of Fredholm operators. Recall that the idea of extending the index beyond the class of Fredholm operators began already in \cite{P7} and was continued in \cite{P48} and \cite{P49}.

Recently in \cite{AZOUZ}, the authors studied operators   $ T \in L(X)$ such that there exists  $ (M, N)$  a GKD  associated to $T$ satisying $T_{\mid_M}$ is an upper semi-Fredholm (resp. a lower semi-Fredholm, a Fredholm) operator and $T_{\mid_N}$ is quasi-nilpotent, called  upper pseudo semi-B-Fredholm (resp.  lower pseudo semi-B-Fredholm,  pseudo B-Fredholm) operators. The class of pseudo-B-Fredholm  operators studied in \cite{AZOUZ} is in fact a particular case of the class of pseud-B-Fredholm studied  here. However, the classes of  upper (resp. lower) pseudo semi-B-Fredholm studied in \cite{AZOUZ} are different  from the classes of  left (resp. right) pseudo-semi-B-Fredholm studied here. In \cite{AZOUZ},  the authors associated  to the class of operators they studied a nullity, deficiency, ascent and descent. However, their definitions of such quantities are incorrect  and induces a misunderstanding in the formulation of their results. More details will be given in Remark \ref{azouz},

%\textbf{Notations:}
\bremark \label{not1} Unless mentioned otherwise, in all this paper,  we will use the following notations.   For  $T \in L(X)$ and  $ P \in L(X)$ an idempotent element of $L(X),$ we will set  $X_1= R(P), X_2= N(P)= R(I-P),$   $T_1= (PTP)_{\mid_{X_1}}$ and  $T_2= (I-P)T(I-P)_{\mid_{X_2}}.$

%\begin{enumerate}
  $\bullet$  If $P$ commutes with $T,$ then  we have  $T= T_1 \oplus T_2,$  here $T_1= T_{\mid_{X_1}}$ and  $T_2= T_{\mid_{X_2}}.$ In this case $(X_1, X_2) \in Red(T).$

   $\bullet$ If $\pi(P)$ commutes with $\pi(T)$ in  $\mathcal{C}_0(X),$   then we have $T= TP + T(I-P)= PTP + (I-P)T(I-P) + F,$ where $F \in L(X)$ is a finite rank operator.
   So  $T= T_1 \oplus T_2 + F.$ In this case $(X_1, X_2) \in Red(PTP)$ and $(X_1, X_2) \in Red((I-P)T(I-P)).$

 $\bullet$ If $\Pi(P)$ commutes with $\Pi(T)$ in  $\mathcal{C}(X),$   then we have $T= TP + T(I-P)= PTP + (I-P)T(I-P) + K,$ where $K \in L(X)$ is a compact operator.
   So $T= T_1 \oplus T_2 + K.$ In this case $(X_1, X_2) \in Red(PTP)$ and $(X_1, X_2) \in Red((I-P)T(I-P)).$

   $\bullet$ $ I_1= P_{\mid X_1},   I_2= (I-P)_{\mid X_2}.  $

   $\bullet$  When we use the homomorphisms $\pi$ or $ \Pi,$  we  mean the natural projections from a Banach subspace $ E$ (appearing in the context of the use of $\pi $ or $\Pi$ ) of $X$   onto the algebra $\mathcal{C}_0(E)$ or $\mathcal{C}(E).$

%\end{enumerate}
\eremark

\section{ P-Fredholm Operators through $ \mathcal{C}_0(X)$}

In a first step, we will extend some  results known  in the case of the Calkin algebra $\mathcal{C}(X)$  to the case of the algebra $ \mathcal{C}_0(X).$
We begin by proving a similar result to \cite [Lemma 1]{BAR}.

\bet \label{idempotent}
Let $p_0$ be an idempotent element of the algebra $ \mathcal{C}_0(X).$  Then there exists
an idempotent $P \in L(X)$ such that $ \pi(P)= p_0.$
\eet

\bp  Assume that $p \in \mathcal{C}_0(X)$ is an idempotent and let $T\in L(X)$ such that $ \pi(T)= p_0.$ As $ p_0^2=p_0,$ then $ T^2-T$ is an operator of finite rank and      thus the spectrum  $\sigma(T)$ of  $T$ is finite.
If $p_0=0,$  we take $P=0.$ So assume that $p_0\neq 0$  and   let $\sigma_1 = \sigma(T) \setminus  \{ 0\} $ and  $\sigma_2 = \{ 0\}.$ Consider two open disjoints sets  $O_1$ and $O_2$ of $\mathbb{C} $ which covers $\sigma_1 $ and $\sigma_2$   respectively. We assume that $ 1 \notin O_2.$
Let $f_1$ and $f_2$ be the functions defined on  $O= O_1 \cup O_2$   by

$ f_1(\lambda) = \left\{
    \begin{array}{ll}
        1 & \mbox{if }  \lambda  \in O_1 \\
        0 &  \mbox{if } \lambda  \in O_2
    \end{array}
\right.$
$\hspace{0.1mm}, \hspace{2cm}  f_2(\lambda) = \left\{
    \begin{array}{ll}
        0 & \mbox{if }  \lambda  \in O_1 \\
        1 &  \mbox{if }  \lambda  \in O_2
    \end{array}
\right.
$ \vspace{2mm}

Let $ P= f_1(T)= \frac{1}{2\Pi i} \int_{\gamma_1} (\lambda I - T) ^{-1}d\lambda.  $  Then $P$ is an idempotent and $ I-P= f_2(T)= \frac{1}{2\Pi i} \int_{\gamma_2}  (\lambda I - T) ^{-1} d\lambda,  $  where  $\gamma_1$  and $\gamma_2$ are respectively simple, closed integration paths,
oriented counterclockwise, which lies in $\rho(T)\cap{O_1} $   and  $\rho(T)\cap{O_2}$  respectively and contain in their interior $\sigma_1$ and   $\sigma_2$ respectively,  $\rho(T)$ being the resolvent set of $T.$

Let
 $ f(\lambda) = \left\{
    \begin{array}{ll}
        \frac{1}{\lambda} & \mbox{if }  \lambda  \in O_1 \\
        0 & \mbox{if } \lambda  \in O_2
    \end{array}
\right.$
\hspace{0.1mm}, then   $f_1(\lambda)= \lambda f( \lambda)$  and so $ P= Tf(T).$

Similarly let $ g(\lambda) = \left\{
    \begin{array}{ll}
        0 & \mbox{if }  \lambda  \in O_1 \\
        \frac{1}{1-\lambda} & \mbox{if } \lambda  \in O_2
    \end{array}
\right.$
\hspace{0.1mm}, then   $f_2(\lambda)= (1-\lambda) g( \lambda)$  and so $ I-P= (I-T)g(T).$

 Thus $P= T R=RT$ and  $I-P= (I-T)S=S(I-T),$  where $ R= f(T)$ and $S= g(T).$
  Hence $ T-P= T(I-P) + (T-I)P= ( T^2- T) ( R- S)$  is a finite rank operator and so $ \pi(P)= p_0.$

We mention that all the functions $f_1, f_2, f, g$ are holomorphic on $O.$ \ep

\bet \label{left invertibility}Let $ T \in L(X).$ Then the following properties are equivalent:
\begin{enumerate}

\item $\pi(T)$ is left invertible in $\mathcal{C}_0(X)$

\item $\Pi(T)$ is left invertible in $ \mathcal{C}(X)$

\end{enumerate}

\eet
\bp
$1) \Rightarrow 2)$ Since every finite rank operator is compact, this implication is trivial.

$2) \Rightarrow 1)$ Suppose  that   $\Pi(T)$ is left invertible in $ \mathcal{C}(X).$  From \cite[Theorem 4.3.3]{CPY},
 there exists a continuous projection $P$ from $X$ onto $R(T).$ Since $N(T)$ is finite dimensional,  $X = N(T) \oplus E$, where $E$ is a closed subspace of $X.$
As $T(E) = R(T), T(E)$  is closed and if  $ T_E: E \rightarrow R(T) $  is the restriction of $T$ to $E$ onto $R(T),$  then $T_E$ is a continuous invertible  operator from $E$  onto $R(T).$  Consider now the operator $U= (T_E)^{-1}PT,$
then $N(U)= N(T)$ is finite dimensional and $ R(U)= E$ is of finite codimension. Therefore $U$ is a Fredholm operator and by the Atkinson theorem \cite[Theorem 1.53]{AI} , there exists  $ S \in L(X),  F_0 \in F_0(X),$
such that $ SU= I + F_0.$ Thus  $S(T_E)^{-1}PT= I+F_0$ and  $\pi(T)$ is left  invertible in $\mathcal{C}_0(X).$
\ep

\bet \label{right invertibility} Let $ T \in L(X).$ Then the following properties are equivalent:
\begin{enumerate}

\item $\pi(T)$ is right invertible in $\mathcal{C}_0(X).$

\item $\Pi(T)$ is right  invertible in $ \mathcal{C}(X).$

\end{enumerate}

\eet
\bp$1) \Rightarrow 2)$ Since every finite rank operator is compact, this implication is trivial.

$2) \Rightarrow 1)$ Suppose  that   $\Pi(T)$ is right invertible in $ \mathcal{C}(X).$  From \cite[Theorem 4.3.2]{CPY},
 there exists a continuous projection
of $X$ onto $N(T).$ Write $X = N(T) \oplus W$, where $W$ is a closed subspace of $X.$
Since $T(W) = R(T), T(W)$  is closed and if  $ T_W: W \rightarrow R(T) $  is the restriction of $T$ to $W$ onto $R(T),$  then $T_W$ is a continuous invertible  map from
$W$  onto $R(T).$  Since $ T \in \Phi_r(X),$   then $R(T) \oplus  F = X$  for a finite dimensional
subspace $F$ of $X.$ Therefore, there exists a bounded projection $P$
of $X$ onto $R(T)$  and the operator $ (T_W)^{-1}P \in L(X).$  Consider now the operator $V= T(T_W)^{-1}P,$
then $N(V)= F$ is finite dimensional and $ R(V)= R(T)$ is of finite codimension. Therefore $V$ is a Fredholm operator and by the Atkinson theorem \cite[Theorem 1.53]{AI} , there exists  $ R \in L(X),  F_0 \in F_0(X),$
such that $ VR= I + F_0.$ Thus  $T(T_W)^{-1}P R= I+F_0$ and  $\pi(T)$ is right  invertible in $\mathcal{C}_0(X).$
\ep

Using Theorem \ref{left invertibility} and Theorem \ref{right invertibility}, we obtain immediately the following result.

\bcor \label{cor-left inv}
\begin{enumerate} Let $T \in L(X).$ Then

\item $T$ is a  left  semi-Fredholm operator if and only if $\pi(T)$ is left invertible in $ \mathcal{C}_0(X).$

 \item $T$ is a  right  semi-Fredholm operators if and only if  $\pi(T)$ is right invertible in $ \mathcal{C}_0(X).$

\end{enumerate}
\ecor

\bremark

  If $p_0$ is an idempotent element of the algebra  $\mathcal{C}_0(X), $ then from Theorem \ref{idempotent}, there exists an idempotent $P \in L(X)$ such that $ \pi(P)=p_0.$ Let $ p= \Pi(P),$    then $p$ is an idempotent of  $\mathcal{C}(X).$
  Conversely if  $p$ is an idempotent of $\mathcal{C}(X),$ then  from \cite[Lemma 1]{BAR},  there exists an idempotent $P \in L(X)$ such that $ \Pi(P)= p.$ Let $ p_0= \pi(P),$    then $ p_0$ is an idempotent of  $\mathcal{C}_0(X).$ In both  cases $ p_0 \cap p \neq \emptyset $ since $ P \in p_0 \cap p. $ (Here $p$ and $p_0$ as equivalence classes can be considered as  subsets of $ L(X)).$

Conversely if   an idempotent $p_0 \in \mathcal{C}_0(X) $ and  an idempotent $p \in \mathcal{C}(X) $ satisfies $ p_0 \cap p \neq \emptyset,$  then there exists  $S \in L(X),$ such that   $p_0= \pi(S)$ and $ p= \Pi(S).$ From Theorem \ref{idempotent}, there exists an idempotent $P \in L(X)$ such that $ \pi(P)= p_0.$  Then  $\Pi(P)= \Pi(S)= p.$ Thus the idempotent $p_0$ and the idempotent $p$  are linked by the idempotent $P$ of $L(X).$

In the sequel, if $ p_0 \cap p \neq \emptyset,$   $P$ will always designate an idempotent of $L(X)$ linking the idempotent $p_0 \in \mathcal{C}_0(X)$ and the idempotent  $ p \in \mathcal{C}(X), $ that's $ \pi(P)= p_0$ and $ \Pi(P)= p.$

\eremark
Recall the following definition from \cite{HB}.

\bdf \label{P-Fredholm}  Let $ T \in L(X)$   and let $p= \Pi(P)$ be an idempotent element of  $\mathcal{C}(X).$
We will say that:

\begin{enumerate}
\item $T$ is a left semi-$P$-Fredholm operator if $ \Pi(T)  $ is left $p$-invertible in $\mathcal{C}(X),$

\item $T$ is a right  semi-$P$-Fredholm operator if $ \Pi(T)$ is right $p$-invertible in $\mathcal{C}(X),$

\item  $T$ is a semi-$P$-Fredholm operator if $T$ is left or right  semi-$P$-Fredholm,

\item  $T$ is a $P$-Fredholm operator if $ \Pi(T)$ is  $p$-invertible in $\mathcal{C}(X).$

\end{enumerate}

\edf

\noindent The classes of left (resp. right) semi-$P$-Fredholm  and $P$-Fredholm operators in Definition\ref{P-Fredholm},   are defined via the Calkin algebra  $\mathcal{C}(X).$ In the next theorem, we will compare these classes to the corresponding classes defined via the algebra  $\mathcal{C}_0(X).$
 Recall that for $  S, T \in L(X),$ the commutator $[S, T]=ST-TS.$  So using  Theorem \ref{right  invertibility} and  Theorem~\ref{left invertibility}, we have

\bet \label{left-right} Let $ P \in L(X)$    be an idempotent element of $L(X).$ Then
\begin{enumerate}

\item   $ \pi(T)  $ is left $p_0$-invertible in $\mathcal{C}_0(X)$  if and only if   $T$ is a left semi-$P$-Fredholm operator and the commutator $[T,P]$ is of  finite rank.

\item  $ \pi(T)  $ is right $p_0$-invertible in $\mathcal{C}_0(X)$  if and only if   $T$ is a right semi-$P$-Fredholm operator and the commutator $[T,P]$ is of  finite rank.

\item   $ \pi(T)  $ is  $p_0$-invertible in $\mathcal{C}_0(X)$  if and only if   $T$ is a $P$-Fredholm operator and the commutator $[T,P]$ is of  finite rank..

\end{enumerate}
\eet

\bp   Let us prove $1).$  Assume that $ \pi(T)  $ is left $p_0$-invertible in $\mathcal{C}_0(X).$  Then the commutator  $[T, P]$ is of finite rank and so it is a compact operator. Moreover  $\pi(T+P)$ is left invertible in  $\mathcal{C}_0(X)$ and its left inverse $ \pi(S)$ commutes with $\pi(T).$ So   $\Pi(T+P)$ is left invertible in  $\mathcal{C}(X)$ and its left inverse $ \Pi(S)$ commutes with $\Pi(T).$  Hence $T$ is a left semi-$P$-Fredholm operator.

Conversely assume  that $T$ is a left semi-$P$-Fredholm operator and the commutator $ [T, P]$   is of  finite rank.
  Then

%\begin{itemize}
  $\bullet$  $\Pi(P) \Pi(T) = \Pi(T) \Pi(P)$

   $\bullet$  There exists  $S \in L(X)$  such that $\Pi(P) \Pi(S) = \Pi(S) \Pi(P)$ and  $\Pi(S) ( \Pi(T) + \Pi(P)) = \Pi(I).$

%\end{itemize}

\noindent As $\Pi(S) ( \Pi(T) + \Pi(P)) = \Pi(I),$ then  $ \Pi(S_1) (\Pi(T_1+ I_1) = \Pi(I_1)$ and
$\Pi(S_2)\Pi(T_2)= \Pi(I_2).$  Hence
$T_1+ I_1$  and    $T_2$ are left semi-Fredholm operators.
From Theorem  \ref{left invertibility}, it follows that $\pi(T_1 + I_1)$ is left invertible in  $\mathcal{C}_0(X_1)$, the  algebra whose unity is $\pi(I_1).$  Similarly     $\pi(T_2)$ is  left invertible in  $\mathcal{C}_0(X_2),$ the  algebra whose unity is $ \pi(I_2).$ Then there exists $R_1 \in L(X_1)$ and  a finite rank operator  $F_1 \in L(X_1)$  such that $R_1( T_1 + I_1)= F_1 +I_1 .$ As the commutator $ [R_1, I_1]=0,$ it is of finite rank. Similarly,  there exists there exists $R_2 \in L(X_2)$ and a finite rank operator   $F_2 \in L(X_2)$  such that $R_2T_2= F_2 +I_2 $ and  $ [R_2, I_2]=0 $ is of finite rank.

 Since $\pi(T)$ and $\pi(P)$ commutes, we have that $\pi(T)+ \pi(P)= \pi(T_1 + I_1) + \pi(T_2),$ because $P_{\mid X_2}=0.$  Hence $\pi(T)+ \pi(P) $ is left invertible in  $\mathcal{C}_0(X), $  having $\pi(R_1 \oplus R_2) $ as a left inverse. Moreover $\pi(R_1 \oplus R_2) $ commutes with $ p_0= \pi(P)$ and we know already that $\pi(T)$ commutes with $ p_0= \pi(P).$ Thus $\pi(T)$ is left $p_0-$ invertible $\mathcal{C}_0(X).$

 The equivalences $2)$ and $3)$ can be proved similarly. \ep

\section{  Pseudo B-Fredholm Operators through $ \mathcal{C}_0(X).$ }

\bdf  \cite[Definition 2]{HA}  An element $a$  of a   ring A with
a unit $e$ is quasinilpotent if, for every x commuting with a, $e -xa $ is invertible in $A.$ The set of all quasinilpotent
elements in a   ring $A$ will be denoted by $QN(A).$

%The set of all quasinilpotent
%elements in a   ring $A$ will be denoted by $QN(A).$

\edf

\bet \label{quasinilpotent} Let $T \in L(X).$ Then the following conditions are equivalent:

\begin{enumerate}

\item  $\Pi(T) $   is quasinilpotent in $ \mathcal{C}(X).$

\item   $\pi(T)$ is quasinilpotent in $ \mathcal{C}_0(X).$

\item    For all   scalar $\lambda \in \mathbb{C}, $     $\pi(I)-\lambda \pi(T) $ is invertible in the algebra $ \mathcal{C}_0(X).$

\end{enumerate}

\eet

\bp
$ 1) \Rightarrow 2) $ Assume that  $\Pi(T) $   is quasinilpotent in $ \mathcal{C}(X)$ and let $ x= \pi(S)$ be an element of $ \mathcal{C}_0(X)$  commuting with $ \pi(T).$ Thus $ TS-ST$ is a finite rank operator, so it is compact and $\Pi(S)$  commutes with $ \Pi(T).$ As  $\Pi(T) $   is quasinilpotent in $ \mathcal{C}(X)$, then $\Pi(S)\Pi(T)$ is also quasinilpotent in $ \mathcal{C}(X).$ Therefore $ I-ST$  is a Fredholm operator and then $ \pi(I)- \pi(S)\pi(T)$ is invertible in $ \mathcal{C}_0(X).$ So $\pi(T)$ is quasinilpotent in $ \mathcal{C}_0(X).$

$ 2) \Rightarrow 3)$ As $\pi(T)$ is quasinilpotent in $ \mathcal{C}_0(X),$  then for all  for all  $x \in \mathcal{C}_0(X) $ commuting with $\pi(T),$ we
have $\pi(I) -x \pi(T)$   is invertible in $\mathcal{C}_0(X).$ This true in particular for $x= \lambda \pi(I).$

$ 3) \Rightarrow 1)$  Assume that for all   scalar $\lambda, $   the element  $\pi(I)-\lambda \pi(T) $ is invertible in the algebra $ \mathcal{C}_0(X).$ Then $\Pi(I)-\lambda \Pi(T) $ is also invertible in the algebra $ \mathcal{C}(X)$ and the spectrum of $\Pi(T)$  in   $\mathcal{C}(X)$ is reduced to $\{0\}.$  Thus $ \Pi(T)$ is quasinilpotent in $\mathcal{C}(X).$   \ep
Recall that in a ring $A$ with a unit and $a \in A,$      $comm(a)= \{ x \in A \mid xa= ax \}$   and $ comm^2(a)=\{ x \in A \mid xy= yx$ for all  $y \in comm(a)\}.$ 

In the next definition, following \cite[Definition 2.2]{KP} where the notion of generalized Drazin invertibility was  already defined in rings , we introduce now the notion of left and right generalized Drazin invertibility in $\mathcal{C}_0(X)$  and  $\mathcal{C}(X).$

\bdf \label{gdrazin}   Let  $T \in L(X).$
\begin{enumerate}

\item  We will say that $\pi(T)$ is generalized Drazin invertible (resp. left generalized Drazin invertible, resp. right  generalized Drazin invertible) in $\mathcal{C}_0(X),$ if there exists an idempotent $P \in L(X) $   such that  $\pi(P)\in comm^2(\pi(T)),$   $ \pi(P)\pi(T)$ is quasi-nilpotent  and $\pi(T)+ \pi(P)$ is  invertible (resp. $\pi(T)$  left $\pi(P)-$invertible, resp. $\pi(T)$ right $\pi(P)-$ invertible) in $\mathcal{C}_0(X).$ In each case, we will say that $P$ is  associated to $T$  in $\mathcal{C}_0(X).$

\item We will say that $\Pi(T)$ is generalized Drazin invertible (resp. left generalized Drazin invertible, resp. right  generalized Drazin invertible) in $\mathcal{C}(X),$ if there exists an idempotent $P \in L(X) $   such that  $\Pi(P)\in comm(\Pi(T)),$   $ \Pi(P)\Pi(T)$ is quasi-nilpotent  and $\Pi(T)+ \Pi(T) $ is  invertible (resp. $\Pi(T)$  left $\Pi(P)-$invertible, resp. $\Pi(T)$ right $\Pi(P)-$ invertible) in $\mathcal{C}(X).$ In each case, we will say that $P$ is  associated to $T$  in $\mathcal{C}(X).$

\end{enumerate}

\edf

\bremark

\begin{enumerate}

\item  It follows from \cite[Theorem 4.2]{KP} that the conditions of the Definition \ref{gdrazin} for the generalized Drazin invertibility in the algebra   $\mathcal{C}_0(X)$  of  $ \pi(T), T \in L(X),$ are equivalent to the existence of
   an element $\pi(S) \in \mathcal{C}_0(X)$  such that:  $$ \pi(S) \in comm^2(\pi(T)),
\pi(S)=\pi(T)\pi(S)^2, \pi(T)-\pi(T)^2\pi(S) \in QN(A).$$

The condition  $ \pi(S) \in comm^2(\pi(T)$ is necessary for the uniqueness of the genaralized Drazin inverse, see \cite[Remark 4.8]{KP}.

 \item   The fact that $\mathcal{C}(X) $ is a  Banach algebra justifies  the relaxed condition $\Pi(T)\in comm(\Pi(T))$  in   Definition \ref{gdrazin} \,  for the generalized Drazin invertibility    of  $ \Pi(T), T \in L(X),$ see  \cite[Remark 4.8]{KP}   for more details. 
\end{enumerate}
\eremark

\bdf

\begin{enumerate}
  \item $T \in L(X)$  is called a Riesz \textbf{-} left semi-Fredholm operator if
 there exists $(M, N) \in Red(T)$ such that $
T_{\mid M}$ is a Riesz   operator and $T_{\mid N} \in \Phi_l(N) .$
  \item  $T \in L(X)$  is called a Riesz \textbf{-} right  semi-Fredholm operator if
 there exists $(M, N) \in Red(T)$ such that $
T_{\mid M}$ is a Riesz  operator and $T_{\mid N} \in \Phi_r(N) .$

  \item  $T \in L(X)$  is called a Riesz \textbf{-}Fredholm operator if
 there exists $(M, N) \in Red(T)$ such that $
T_{\mid M}$ is a Riesz  operator and $T_{\mid N}$ is a Fredholm operator.
\end{enumerate}
\edf

\bet \label{first equivalence} Let $T \in L(X),$   then the following properties
are equivalent.

\begin{enumerate}

\item $\pi(T)$ is generalized Drazin invertible in $\mathcal{C}_0(X).$

\item $\pi(T)$  is left and  right generalized Drazin invertible in $\mathcal{C}_0(X).$

\end{enumerate}
\eet

\bp Since the implication $ 1) \Rightarrow 2)$ is trivial, we only need to prove the implication $ 2) \Rightarrow 1).$

So assume that $\pi(T)$  is left and  right generalized Drazin invertible in $\mathcal{C}_0(X).$ Then there exists two idempotents $P , Q  \in L(X)$ such that $\pi(P) \in comm^2(\pi(T))$, $\pi(Q) \in comm^2(\pi(T)),$  $\pi(T)$ is left $\pi(P)-$invertible in $\mathcal{C}_0(X),$ $\pi(T)\pi(P) $ is quasinilpotent in $\mathcal{C}_0(X),$   $\pi(T)$ is right $\pi(Q)-$invertible in $\mathcal{C}_0(X)$ and $\pi(T)\pi(Q) $ is quasinilpotent in $\mathcal{C}_0(X).$ As $\pi(T)$ is left $\pi(P)-$invertible in $\mathcal{C}_0(X),$ it has a left inverse  $\pi(S), S\in L(X)$ such that $\pi(S)\pi(P)= \pi(P) \pi(S).$ Then  $\pi(I-P)\pi(T)\pi(I-P)$ is left invertible in the algebra $\pi(I-P)\mathcal{C}_0(X)\pi(I-P)$ whose identity element is $ \pi(I-P), $  having as left inverse $\pi(I-P) \Pi(S)\pi(I-P)$ which commutes with $\pi(P).$
So $$\pi(I-P) \Pi(S)\pi(I-P) \pi(I-P)\pi(T)\pi(I-P)= \pi(I-P).$$
Since $\pi(T)\pi(Q) $ is quasinilpotent in $\mathcal{C}_0(X)$  and $\pi(P) \in comm^2(\pi(T)),$   then from  Theorem \ref{quasinilpotent}, $\Pi(I-P)\Pi(T)\Pi(Q) $ is quasinilpotent in $\mathcal{C}(X).$  As $\pi(I-P) \pi(S)\pi(I-P)\pi(I-P) \pi(T)\pi(I-P)= \pi(I-P),$ then  $\Pi(I-P) \Pi(S)\Pi(I-P)\Pi(I-P) \Pi(T)\Pi(I-P)= \Pi(I-P).$ Hence  $\Pi(I-P) \Pi(S)\Pi(I-P) \Pi(T)\Pi(Q) \Pi(I-P)= \Pi(I-P)\Pi(Q).$ Thus $\Pi(I-P) \Pi(S^n)\Pi(I-P) \Pi(T^n)\Pi(Q) \Pi(I-P)= \Pi(I-P)\Pi(Q),$  for every positive integer $n,$  because $\Pi(S)$ commutes with $\Pi(P).$ As  $ \Pi(I-P)\Pi(Q)$  is an idempotent, then if $\Pi(I-P)\Pi(Q)\neq 0,$ we have     $1\leq  \| \Pi(I-P)\Pi(Q) \|.$ Therefore
$$1 \leq  \| \Pi(I-P) \Pi(S^n)\Pi(I-P) \Pi(T^n)\Pi(Q) \Pi(I-P)\|  $$
$$ \leq  \| \Pi(I-P)\| \| \Pi(S\|^n \| \Pi(I-P)\| \|\Pi(T^n)\Pi(Q)\Pi(I-P)\|.$$

Hence
$$ 1 \leq \| \Pi(I-P)\|^{ \frac{1}{n}} \| \Pi(S\| \| \Pi(I-P)\|^{ \frac{1}{n}} \|\Pi(T^n)\Pi(Q)\Pi(I-P)\|^{ \frac{1}{n}}.$$

But this impossible since  $\Pi(T)\Pi(Q)\Pi(I-P)$ is quasinilpotent in $\mathcal{C}(X).$ Thus we must have  $\Pi(I-P)\Pi(Q)= 0$  and so $\Pi(Q)= \Pi(PQ).$
Similarly we have $\Pi(P)= \Pi(QP)= \Pi(PQ),$  hence $ \Pi(P)= \Pi(Q).$

Moreover as $\pi(T)+ \pi(Q)$ is right invertible in $\mathcal{C}_0(X),$ then from Theorem \ref{right invertibility},  $\Pi(T)+ \Pi(P)= \Pi(T) +\Pi(Q)$ is right invertible in $\mathcal{C}(X).$  From Theorem \ref{left invertibility}, it follows that $\pi(T)+ \pi(P)$ is right invertible in $\mathcal{C}_0(X).$ As we know already that $\pi(P) \in comm^2(\pi(T)),$  $\pi(T)\pi(P) $ is quasinilpotent in $\mathcal{C}_0(X)$  and $ \pi(T)+\pi(P)$ is left invertible in $\mathcal{C}_0(X),$ then $\pi(T)$ is  generalized Drazin invertible in $\mathcal{C}_0(X).$ \ep

Recall the following definition from \cite{HB}:

\bdf \label{pseudo} Let $ T \in L(X).$    We will say that:

\begin{enumerate}
  \item $T$ is a left pseudo   semi-B-Fredholm operator if $ \Pi(T)  $ is left generalized Drazin  invertible in $\mathcal{C}(X).$
  \item $T$ is a right pseudo semi-B-Fredholm operator if $ \Pi(T)$ is right generalized  Drazin invertible in $\mathcal{C}(X).$
  \item $T$ is a  pseudo semi-B-Fredholm operator if it is right  or left pseudo semi-B-Fredholm.
  \item $T$ is a pseudo B-Fredholm operator if $ \Pi(T)$ is  generalized Drazin invertible in $\mathcal{C}(X).$
\end{enumerate}

\edf

\begin{remark}

\begin{enumerate}

\item The class of  pseudo B-Fredholm  was defined in  \cite[Definition 2.4]{P45}. It involves the class
of B-Fredholm operators defined in \cite{P7}.

 \item In \cite{BO}, in 2015, the author studied  generalized Drazin invertible elements under Banach algebra homomorphisms.  In the case of the  Calkin Algebra,  he proved in \cite[Theorem 6.1]{BO} that an operator whose image is generalized Drazin invertible in the Calkin algebra is a Riesz-Fredholm operator. Later in 2019, in \cite[Theorem 2.10]{P45},  the authors proved the same result following a shortened proof. 

 Here in this paper,  the concept of p-invertibility enable us to give left and right version of \cite[Theorem 6.1]{BO}.
 
\end{enumerate}

\end{remark}

\noindent The classes of left (resp. right) pseudo semi-B-Fredholm  and pseudo $B$-Fredholm operators in Definition \ref{P-Fredholm},   are defined via the Calkin algebra  $\mathcal{C}(X).$ In the next theorems, we will compare these classes to the corresponding classes defined via the algebra  $\mathcal{C}_0(X).$

\bet \label{left-pseudo-0} Let $T \in L(X).$ Then the following properties
are equivalent:

\begin{enumerate}

\item $\pi(T)$ is left generalized Drazin invertible  in $ \mathcal{C}_0(X).$
\item  $T$ is left pseudo   semi-B-Fredholm operator and there exists an idempotent $P$ associated to $\Pi(T)$ in $ \mathcal{C}(X)$  such that $\pi(P) \in comm^2( \pi(T)).$
\item  There exists an idempotent $P \in L(X)$ such that $\pi(P) \in comm^2( \pi(T)),$ $T_1 $ is a Riesz operator and  $T_2 $ is a semi-Fredholm operator.

In this case $T$ is a finite rank perturbation of a Riesz-left semi-Fredholm operator.
\end{enumerate}

\eet
\bp  $1)\Rightarrow 2)$  Assume that  $\pi(T)$ is  left generalized Drazin invertible in  $\mathcal{C}_0(X).$ Then there exist an
idempotent $P \in L(X)$   associated to $ T$ in $\mathcal{C}_0(X),$  so:

%\begin{itemize}
  $\bullet$  $p_0 \in comm^2(\pi(T)) $, in particular $p_0 \pi(T) = \pi(T) p_0.$

   $\bullet$ $p_0 \pi(T)$ is quasinilpotent in $\mathcal{C}_0(X).$

   $\bullet$  There exists  $S \in L(X)$  such that $p_0 \pi(S) = \pi(S) p_0$ and  $\pi(S) ( \pi(T) + p_0)) = \pi(I).$

%\end{itemize}

Let $ p= \Pi(P),$ as $p_0 \in comm^2(\pi(T))$ and $ p_0 \pi(T)$ is quasinilpotent in $\mathcal{C}_0(X),$ then $ p \Pi(T) = \Pi(T) p$
 and   from Theorem \ref{quasinilpotent},  $ p \Pi(T)$ is quasinilpotent in $\mathcal{C}(X).$
 Moreover we have   $p \Pi(S) = \Pi(S) p$ and  $\Pi(S) ( \Pi(T) + p)) = \Pi(I).$ Therefore $ \Pi(T)$ is left generalized Drazin invertible in $\mathcal{C}(X) $ and so $T$ is left pseudo   semi-B-Fredholm operator, the  idempotent $P$ is  associated to $\Pi(T)$ in $ \mathcal{C}(X)$ and $\pi(P) \in comm^2( \pi(T)).$

 $2)\Rightarrow 3)$  Assume that   $T$ is left pseudo semi-B-Fredholm operator and there exists an idempotent $P$ associated to $T$ in $ \mathcal{C}(X)$  such that $\pi(P) \in comm^2( \pi(T)).$
 Since $\pi(T)$ and $\pi(P)$ commutes, we have  $\pi(PTP)=\pi(TP)$ and  from Theorem \ref{quasinilpotent}, it follows that $PTP$
 is a Riesz operator.  As $\pi(T)$ and $\pi(P)$ commutes, from Remark \ref{not1}, we have:

 \begin{equation}\label{Decomposition}
 T=T_1\oplus T_2 +F,
 \end{equation}

   where $F$ is a finite rank operator.

 As $PTP$  is a Riesz operator, then $T_1$ is also Riesz operator. Moreover  we show that $T_2$ is a left semi-Fredholm operator.
Since $P$ is associated to $T$ in $ \mathcal{C}(X),$  then $\Pi(T)$  is left $\Pi(P)$-invertible in $ \mathcal{C}(X).$ From Theorem \ref{left invertibility}, it follows that $\pi(T) +\pi(P)$ is left invertible in $ \mathcal{C}_0(X).$   So  there exists  $S \in L(X)$  such that   $\pi(S) ( \pi(T) + \pi(P))) = \pi(I).$ As $\pi(T) \pi(I-P)=\pi(I-P)\pi(T),$  then   we have: $ (I-P)S(I-P) (I-P) (T+P)(I-P)= I-P  + (I-P)F'(I-P)$  and $ (I-P)S(I-P)_{\mid{X_2}}  T_2= I_2  + (I-P)F'(I-P)_{\mid X_2},$
where   $F'$ is a  finite rank operator.
Hence $T_2$ is a left semi-Fredholm operator. According to (\ref{Decomposition}), we see that  $T$ is a finite rank perturbation of
a Riesz-left  semi-Fredholm operator.

$ 3)\Rightarrow 1)$ Conversely assume that there exists an idempotent $P \in L(X)$ such that $\pi(P) \in comm^2( \pi(T)),$   $T_1 $ is a Riesz operator and  $T_2$ is a left semi-Fredholm operator. We have $ T= T_1 \oplus T_2 +F,$ where $F$ is a finite rank operator.
So  $ T_1\oplus T_2 +P = ( T_1  \oplus 0) +  P  + ( 0 \oplus T_2) = P[ (T_1 \oplus 0)+ I)] P  + (I-P) ( I_1 \oplus  T_2 ) (I-P).$  As $T_1$ is a Riesz operator, then $T_1 \oplus 0$ is also a Riesz operator and $\pi( (T_1 \oplus 0)+ I))=  \pi( (T_1 \oplus 0)) + \pi(I)$ is invertible in $\mathcal{C}_0(X).$ Let $ \pi(S_1) $ be its inverse, where $ S_1 \in L(X).$

As $T_2$ is a left  semi-Fredholm operator in $L(X_2),$  then from Corollary \ref{cor-left inv}, there exists $ S_2 \in  L(X_2)$  such that $ S_2 T_2- I_2$ is a finite rank operator. Moreover $I_1 \oplus S_2$ commutes with $P$   because
$(I_1  \oplus S_2)P=P(I_1 \oplus S_2)=I_1 \oplus 0.$
We observe that $\pi( I_1\oplus  T_2) $ is left invertible in  $\mathcal{C}(X)$ having  $ \pi( I_1 \oplus S_2)$ as a left inverse. Then:
$$ \pi((PS_1P + (I-P)(I_1\oplus S_2) (I-P))) \pi ( T+P)$$ $$=\pi((PS_1P + (I-P)(I_1\oplus S_2) (I-P))) \pi (T_1\oplus T_2 +P)$$
$$=\pi([PS_1P + (I-P)(I_1\oplus S_2) (I-P)]) \pi ( P[ (T_1 \oplus 0)+ I)] P  + (I-P) ( I_1 \oplus  T_2 ) (I-P))=\pi(I).$$
It is easily seen that  $\pi((PS_1P + (I-P)( I_1 \oplus S_2) (I-P)))$ commutes with $\pi(P).$

\noindent Moreover we have $\pi(PT)= \pi( T_1 \oplus 0)$   is  quasinilpotent in $\mathcal{C}_0(X)$  because $ T_1 \oplus 0$ is a Riesz operator.  Thus $ T$ is a left generalized Drazin invertible in $\mathcal{C}_0(X).$
\ep

%\noindent As $T =  T_1\oplus T_2 +F,$ then $T$ is a finite rank perturbation of a Riesz-left semi-Fredholm operator.

 The next result is given without proof, due to the similarity of its proof with that of Theorem \ref{left-pseudo-0}.

\bet \label{right-pseudo-0} Let $T \in L(X).$ Then the following properties
are equivalent:

\begin{enumerate}

\item $\pi(T)$ is right generalized Drazin invertible  in $ \mathcal{C}_0(X).$

\item  $T$ is right pseudo   semi-B-Fredholm operator and there exists an idempotent $P$ associated to $\Pi(T)$ such that $\pi(P) \in comm^2( \pi(T)).$

\item    There exists an idempotent $P \in L(X)$ such that $\pi(P) \in comm^2( \pi(T)),$ $T_1 $ is a Riesz operator and  $T_2 $ is a right semi-Fredholm operator.

In this case $T$ is a finite rank perturbation of a Riesz-right semi-Fredholm operator.

\end{enumerate}

\eet

\bet \label{left-pseudo} Let $T \in L(X).$ Then the following properties
are equivalent:

\begin{enumerate}

\item $\pi(T)$ is  generalized Drazin invertible  in $ \mathcal{C}_0(X).$

\item  $T$ is  pseudo B-Fredholm operator and there exists an idempotent $P$ associated to $T$ in $\mathcal{C}(X)$ such that $\pi(P) \in comm^2( \pi(T)).$

\item   There exists an idempotent $P \in L(X)$ such that $\pi(P) \in comm^2( \pi(T)),$ $T_1 $ is a Riesz operator and  $T_2 $ is a Fredholm operator.

In this case $T$ is a finite rank perturbation of a Riesz-Fredholm operator.

\end{enumerate}

\eet

\bp

$1)\Rightarrow 2)$  Assume that  $\pi(T)$ is   generalized Drazin invertible in  $\mathcal{C}_0(X).$ Then there exists an idempotent $P \in  L(X)$ associated to $ T$ in $\mathcal{C}_0(X),$  so:

%\begin{itemize}
  $\bullet$  $p_0 \in comm^2(\pi(T)) $, in particular $p_0 \pi(T) = \pi(T) p_0.$

   $\bullet$ $p_0 \pi(T)$ is quasinilpotent in $\mathcal{C}_0(X).$

   $\bullet$  There exists  $S \in L(X)$  such that %$p_0 \pi(S) = \pi(S) p_0$ and
    $\pi(S) ( \pi(T) + p_0)) = \pi(S) ( \pi(T) + p_0))= \pi(I).$

%\end{itemize}

Let $ p= \Pi(P),$ as $p_0 \in comm^2(\pi(T))$ and $ p_0 \pi(T)$ is quasinilpotent in $\mathcal{C}_0(X),$ then $ p \Pi(T) = \Pi(T) p$
 and   from Theorem \ref{quasinilpotent},  $ p \Pi(T)$ is quasinilpotent in $\mathcal{C}(X).$
 Moreover    $\Pi(S) ( \Pi(T) + p)) =  ( \Pi(T) + p))\Pi(S)=\Pi(I).$ Therefore   $T$ is  pseudo B-Fredholm operator, the  idempotent $P$ is  associated to $\Pi(T)$ in $ \mathcal{C}(X)$ and $\pi(P) \in comm^2( \pi(T)).$

$ 2)\Rightarrow 3)$  Assume that   $T$ is a pseudo B-Fredholm operator and there exists an idempotent $P$ associated to $\Pi(T)$ in $ \mathcal{C}(X)$  such that $\pi(P) \in comm^2( \pi(T)).$
  Then, there exists a finite rank operator $ F \in L(X),$ such that:

\begin{equation}\label{first-2}
    T=T_1\oplus T_2+F.
  \end{equation}

  It's easily seen that $T_1$ is Riesz operator. Let us show that $T_2$ is a Fredholm operator.

 As $P$ is associated to $T$ in $ \mathcal{C}(X),$  then $\Pi(T) +\Pi(P)$ is  invertible in $ \mathcal{C}(X)$ and so  $\pi(T) +\pi(P)$ is  invertible in $ \mathcal{C}_0(X).$   Hence  there exists  $S \in L(X)$  such that   $\pi(S) ( \pi(T) + \pi(P))) =  ( \pi(T) + \pi(P)))\pi(S)= \pi(I).$ Then $\pi(S)\pi(P)=\pi(P)))\pi(S)$ and
 $$ S_2 T_2= I_2 + F_1,\,\,T_2 S_2 = I_2  + F_2,$$
where   $F_1, F_2$ are  finite rank operators in $L(X_2).$
Therefore $T_2$  is a Fredholm operator.

$ 3)\Rightarrow 1)$  Assume that there exists an idempotent $P \in L(X)$ such that $\pi(P) \in comm^2( \pi(T)),$ $T_1 $ is a Riesz operator and  $T_2 $ is a Fredholm operator. From Theorem \ref{left-pseudo-0} and  Theorem \ref{right-pseudo-0}, it follows that $\pi(T)$ is  left generalized invertible  and right generalized invertible in $ \mathcal{C}_0(X).$ Then from Theorem \ref{first equivalence}, it follows that $\pi(T)$  is generalized invertible in $ \mathcal{C}_0(X).$

\noindent As $\pi(P) \in comm^2( \pi(T)),$ then  $T =  T_1\oplus T_2 +F,$  $F$ being a finite rank operator. So $T$ is a finite rank perturbation of a Riesz-Fredholm operator.  \ep

It follows from the properties of generalized Drazin invertibility in a Banach algebra, that if  $T \in L(X)$  is  a   pseudo B-Fredholm operator and if $P$ is an idempotent associated to $T$ in $\mathcal{C}(X),$  then $ comm(\Pi(T)) \subset comm(\Pi(P)) .$ This a direct consequence of the properties of the functional calculus,  see \cite[Theorem 3.1]{Koliha} for the construction of the idempotent $p=\Pi(P).$

\bremark \label{com-com} Let $T \in L(X)$  be a   pseudo B-Fredholm operator and let $P$ be an idempotent associated to $T$ in $\mathcal{C}(X).$
 Then  $\pi(T)$ is  generalized Drazin invertible  in $ \mathcal{C}_0(X)$ if and only if $ comm(\pi(T) \subset comm(\pi(P)).$ Indeed $\pi(P) \in comm^2( \pi(T))$  is equivalent to   $ comm(\pi(T)) \subset comm(\pi(P))$ and then the remark is consequence of Theorem \ref{left-pseudo}.

\eremark

\bex  Let $P$ be an idempotent element of $L(X)$  and let $K$ be a compact operator. Consider the operator $T= I-P+K.$
Then  $T$ is a pseudo B-Fredholm operator with  $P$ as an  associated idempotent  in   $ \mathcal{C}(X).$ From \cite[Theorem 4.2]{Koliha}, the idempotent $ p=\Pi(P)$ is the unique idempotent in  $ \mathcal{C}(X)$ such that $\Pi(T) p= p \Pi(T),$ $\Pi(T) + p$ is invertible and  $\Pi(T)p$ is quasi-nilpotent.

Now if there exists an idempotent  $Q$ in $L(X)$ such that $\pi(T)$ is  generalized Drazin invertible in  $ \mathcal{C}_0(X)$ and $Q$ is associated to $T$  in $ \mathcal{C}_0(X),$ then  $Q$ is associated to $T$ in $ \mathcal{C}(X).$ So we must have $ \Pi(Q)= \Pi(P).$
In particular  if  $comm(\pi(T)) \subset comm( \pi(K)),$ then  $comm(\pi(T)) \subset comm( \pi(P)).$  From Remark \ref{com-com} it follows that $\pi(T)$ is  generalized Drazin invertible  in $ \mathcal{C}_0(X)$ and $P$ is  associated to $T$  in   $ \mathcal{C}_0(X).$

{\bf \underline{Open Question:}} More generally, what are the conditions on $P$ and $K$  so that  $\pi(T)$ is  generalized Drazin invertible in  $ \mathcal{C}_0(X)?$
\eex

\section{Index of Fredholm type Operators}

In this section,  we will  define   an index for  pseudo-B-Fredholm  operators. As we will see, this new index extends  naturally the usual Fredholm index and satisfies similar properties.

\bdf Let $ T \in L(X)$ be a  pseudo semi-B-Fredholm operator and let $P$ be an idempotent element in  $L(X)$  associated to $T$ in $\mathcal{C}(X).$    We define the index  \textit{\textbf{{ind}}}($T$) of $T$ as the index  of the semi-Fredholm operator $T+P, $  that's  \textit{\textbf{{ind}}}($T$)= ind$(T+P).$

\edf

\bet \label{index}

Let $ T \in L(X)$ be a  pseudo semi-B-Fredholm operator.   Then the index  \textit{ind}($T$) of $T$  is well defined and \textit{\textbf{{ind}}}($T$)=  $ \lim \limits_{\lambda \rightarrow 0} ind ( T- \lambda I).$

\eet

\bp  Assume that $T$ is a left pseudo semi-B-Fredholm and  let $P$ be an idempotent element of $L(X)$ associated to $T$ in $\mathcal{C}(X).$    We have $ T+P= (T_1 +I_1) \oplus T_2 +K,$  where $K$ is a compact operator.
As $T_1$ is a Riesz operator, then $T_1 + I_1$ is a Fredholm operator of index $0.$ Moreover $T_2$ is a left semi-Fredholm operator. Thus $(T_1 +I_1) \oplus T_2 $ is an upper semi-Fredholm operator and $ind (T+P)= ind((T_1 +I_1) \oplus T_2)= ind(T_2) .$

Now let $ \lambda \in \mathbb{C}\setminus \{0\},$ then $ T- \lambda I= (T_1 -\lambda I_1) \oplus (T_2- \lambda I_2) +K.$
As $T_1$ is a Riesz operator,  then $ T_1 -\lambda I_1$ is a Fredholm operator of index $0.$  Moreover if $ \mid \lambda \mid $ is small enough, then $T_2- \lambda I_2$ is also a semi-Fredholm operator.

If $T+P$ is a Fredholm operator, then $T_2$ is  a Fredholm operator and for $ \mid \lambda \mid $ small enough,
$T_2- \lambda I_2$ is  also a Fredholm operator and  by the continuity of the index  $ind(T_2- \lambda I_2)= ind(T_2).$    We have   $ind(T+P)= ind(T_2)=  ind(T_2- \lambda I_2)= ind((T_1  - \lambda I_1) \oplus (T_2- \lambda I_2) $ because $T_1  - \lambda I_1$ is a Fredholm operator of index $0.$
Thus \textit{\textbf{{ind}}}($T$)= ind$(T+P)= ind(T_2)= \lim \limits_{\lambda \rightarrow 0} ind ( T- \lambda I).$

Now if  $T+P$ is not a Fredholm operator, then $\beta ( T_2)= +\infty.$ From \cite[Theorem 4.2.1]{CPY}, it follows that for $\mid \lambda \mid $ small enough, we have $\beta(T_2- \lambda I_2)= +\infty.$ Hence \textit{\textbf{{ind}}}($T$)= ind$(T+P)= \lim \limits_{\lambda \rightarrow 0} ind ( T- \lambda I)= -\infty.$

In both cases, the index \textit{\textbf{{ind}}}($T$)=$\lim \limits_{\lambda \rightarrow 0} ind ( T- \lambda I).$ So it is well defined and it is independent of the choice of the associated idempotent $ p=\Pi(P).$

\noindent The proof is similar in the case where  $T$ is a right pseudo semi-B-Fredholm or pseudo B-Fredholm. \ep

\bremark
 Let $ T \in L(X)$ be a  pseudo semi-B-Fredholm operator and let $P$ be an idempotent element in  $L(X)$  associated to $T$ in $\mathcal{C}(X.)$ It follows from the proof of Theorem \ref{index}, that  the index  \textit{\textbf{{ind}}}($T$)= $ ind(T_2).$
\eremark

\bremark \label{azouz}

\begin{enumerate}
\item Clearly if $T$ is a left (resp. right) semi-Fredholm operator, then  $0= \Pi(0)$ is an idempotent element of $\mathcal{C}(X)$ commuting with $\Pi(T),$
 such that $\Pi(T)$ is left (resp. right)  $0$-invertible in $\mathcal{C}(X)$ and $0=\Pi(T) 0$ is quasinilpotent. Thus   the new index \textit{\textbf{{ind}}} coincides with the usual index for these classes of operators. In particular the new  index \textit{\textbf{{ind}}} coincides with the usual index in the case of Fredholm  operators.

\item In \cite{AZOUZ}, an operator is  called  upper pseudo semi-B-Fredholm (resp.  lower pseudo semi-B-Fredholm,  pseudo B-Fredholm)
operator if there exists $ (M, N)$  a GKD associated to $T$ such that $T_{\mid_M}$ is an upper semi-Fredholm (resp. a lower semi-Fredholm, a Fredholm) operator and $T_{\mid_N}$ is
 quasi-nilpotent. In \cite[Definition 2.15]{AZOUZ},  the following quantities were associated to these operators.
  The nullity, the deficiency, the ascent and the descent of such $T$ defined respectively, by

 $ \alpha(T )= dim N( T_{\mid M}), \beta(T)= codim R(T_{\mid M}),
 p(T ) = inf\{n \in \mathbb{N}:  N( (T_{\mid M})^m) =N( (T_{\mid M})^n) $ for all integer $m \geq n$  \} and $q(T ) = inf\{n \in \mathbb{N}:  R( (T_{\mid M})^m) =R( (T_{\mid M})^n) $ for all integer $m \geq n$  \}.

However these definitions are ambiguous  since the same quantities are already defined by:
 $ \alpha(T )= dim N( T), \beta(T)= codim R(T),
 p(T )= inf\{n \in \mathbb{N}:  N( T^m) =N( T^n )$ for all integer $m \geq n$  \} and $q(T ) = inf\{n \in \mathbb{N}:  R( T^m) =R( T^n) $ for all integer $m \geq n$  \}.\\
For example consider the  operator $ T= I \oplus 0,$ on the Banach space $X\oplus X,$  where $X$ is an infinite dimensional Banach space. Following \cite[Definition 2.15]{AZOUZ} we obtain $\alpha(T )= 0,$  however the nullity of $T$ in the usual sense is equal to $\infty.$  The deficiency of $T$   in the sense of \cite[Definition 2.15]{AZOUZ} is $\beta(T )= 0,$ however the deficiency of $T$ in the usual sense is equal to $\infty.$
The ascent of $T,$   in the sense of \cite[Definition 2.15]{AZOUZ} is $p(T )= 0,$ however the ascent of $T$ in the usual sense is equal to $1.$ Similarly the descent of $T,$   in the sense of \cite[Definition 2.15]{AZOUZ} is $q(T )= 0,$ however the descent of $T$ in the usual sense is equal to $1.$

To make things clear in \cite{AZOUZ},  it enough to consider that the quantities defined in  \cite[Definition 2.15]{AZOUZ} are  relatively to the subspace $M$ as  in Definition \ref{relative-index}.

%as a sample, how one could interpret their  results \cite[Theorem 2.19]{AZOUZ} \cite[Theorem 2.21]{AZOUZ}.

\item For a pseudo semi-B-Fredholm  $T$ as defined  in \cite{AZOUZ}, the authors associated an index by setting $ind(T)= ind(T_{\mid M}).$ This definition is the same as our definition of the index of pseudo semi-B-Fredholm operators, for if $P$ is the projection on $N$ in parallel to $M,$  then $ind(T_{\mid_M})=ind(I-P)T(I-P)_{\mid N(P)}= \alpha_M(T)- \beta_M(T), $ where  $\alpha_M(T)$  and  $\beta_M(T)$ are respectively the nullity and the deficiency of $T$ relatively to $M.$

\end{enumerate}

\eremark

The following result is obvious, however it shades light on the difference between the pseudo semi-B-Fredholm operators studied in \cite{AZOUZ} and the pseudo semi-B-Fredholm studied here.
\bet  Let $ T \in L(X).$ Then the following properties are equivalent:
\begin{enumerate}

\item There exists $ (M, N)$  a GKD associated to $T$ such that $T_{|_M}$ is an upper semi-Fredholm (resp. a lower semi-Fredholm, a Fredholm) operator and $T\mid_N$ is
 quasi-nilpotent.

\item  There exists an idempotent $ P$ in $L(X)$ commuting  with $T,$ such that  $ PT$ is quasi-nilpotent and $(I-P)T$ is an upper semi-Fredholm (resp. a lower semi-Fredholm, a Fredholm) operator.

\end{enumerate}

\eet

In the case of a pseudo semi-B-Fredholm operator $T$, as we defined it here,  if   $ p=\Pi(P)$ is an idempotent an associated to $T,$ then $P$ is not necessarily commuting with $T.$ Moreover for a  pseudo semi-B-Fredholm operator $T$ as defined in \cite{AZOUZ}, $T_{\mid M}$ is an upper semi-Fredholm (resp. a lower semi-Fredholm) but not necessarily a left (resp. right) semi-Fredholm operator.

We give now some results which shows that the new index defined here,  satisfies also several known properties as the usual Fredholm index.

\bprop Let $ T \in L(X)$ be a  pseudo semi-B-Fredholm (resp. pseudo B-Fredholm)  operator and let $K \in L(X)$ be a compact operator.  Then $ T+K$ is a pseudo semi-B-Fredholm (resp. a pseudo B-Fredholm) operator and
\textit{\textbf{{ind}}}($T$)= \textit{\textbf{{ind}}}($T+ K)$).
\eprop
\bp As $\Pi(T+K)= \Pi(T)$ and $T$ is pseudo semi-B-Fredholm (resp. pseudo B-Fredholm)  operator, then $T+K$ is also  a  pseudo semi-B-Fredholm (resp. pseudo B-Fredholm)  operator. 
 We have      \textit{\textbf{{ind}}}($T$)= $\lim \limits_{\lambda \rightarrow 0} ind ( T- \lambda I)= \lim \limits_{\lambda \rightarrow 0} ind ( T+K- \lambda I) $ because $ T-\lambda I$ is a semi-Fredholm (resp. pseudo B-Fredholm) operator and $ K$ is compact. Thus \textit{\textbf{{ind}}}($T$)= \textit{\textbf{{ind}}}($(T+K)$).
 \ep

\bet \label{product} Let $S, T$ be two   left (resp. right) pseudo semi-B-Fredholm operators such that there exists an idempotent $P \in L(X)$ associated to both $S$ and $T$  in $\mathcal{C}(X).$ If  the commutator $[S_1, T_1]$ is compact, then $ST$ is a left (resp. right) pseudo-semi-B-Fredholm operator and \textit{\textbf{{ind}}}($ST$)= \textit{\textbf{{ind}}}($S$) + \textit{\textbf{{ind}}}($T$).
\eet

\bp Let $S,T$ be two   left  pseudo semi-B-Fredholm operators and let   $P$ be a common associated idempotent. Set $ S=  S_1 \oplus S_2 +K_1,$ $ T= T_1 \oplus T_2 +K_2,$ where $K_1$ and $K_2$ are compact operators. Then $ ST=S_1T_1\oplus S_2T_2 +K,$ where $K$ is a compact operator.
As     $\Pi(P)$ commutes with $\Pi(S)$ and $\Pi(T),$  then $\Pi(P)$ commutes with $\Pi(ST),$  thus
$S_1T_1= ( PSP) (PTP)_{\mid X_1}= ( PSTP)_{\mid X_1} + K_3= (ST)_1 + K_3,$ where $K_3 \in L(X_1)$ is a compact operator.  As $S_1$  and $T_1$ are  Riesz operators and  the commutator $[S_1, T_1]$ is compact, then  $(ST)_1=  S_1T_1 -K_3$ is a Riesz operator.\\
$(ST)_1$ is also  Riesz operator.\\
 Moreover $S_2 T_2$ is a left semi-Fredholm operator, because $S_2$  and $T_2$ are both left  semi-Fredholm operators.  As  $\Pi(P)$ commutes with $\Pi(S), \Pi(T)$ and  $\Pi(ST),$ then $(ST)_2= S_2T_2 + K_4,$ where $K_4\in L(X_2) $ is a compact operator. Hence  $(ST)_2$ is a left semi-Fredholm operator and  $ind( (ST)_2)=ind( S_2T_2)=  ind( S_2) + ind( T_2 )=$ \textit{\textbf{{ind}}}($S$) + \textit{\textbf{{ind}}}($T$).\\
 We have  $ ST=S_1T_1\oplus S_2T_2 +K= (ST)_1 \oplus (ST)_2 + (K_3\oplus K_4) +K.$ From \cite[Theorem 3.5]{HB}, it follows that $ST$ is a left pseudo semi-B-Fredholm operator and  \textit{\textbf{{ind}}}($ST$) =\textit{\textbf{{ind}}}($S$) + \textit{\textbf{{ind}}}($T$).\\
 Following   similar steps  and using \cite[Theorem 3.6]{HB}  which characterizes right pseudo semi-B-Fredholm operators, we can prove  the case of the product of two right pseudo-semi-B-Fredholm operators. \ep

Theorem \ref{product} gives  immediately the following result.
\bcor
Let $S, T$ be two    pseudo B-Fredholm operators such that there exists an idempotent $P \in L(X)$ associated to both $S$ and $T$  in $\mathcal{C}(X).$  If $S_1$  and $T_1$ are commuting operators, then $ST$ is a pseudo-B-Fredholm operator and \textit{\textbf{{ind}}}($ST$)= \textit{\textbf{{ind}}}($S$) + \textit{\textbf{{ind}}}($T$).
\ecor

\bet \label{left-perturbation} Let $T$ be a   left  pseudo-semi-B-Fredholm operator and let  $P \in L(X)$  be an idempotent associated to $T$ in $\mathcal{C}(X).$  Then there exists an $ \epsilon > 0,$ such that if $U \in L(X)$  is such that:

\begin{enumerate}

\item  $\Pi(P)$ commutes with $ \Pi(U),$
\item  $U_1$ is a Riesz operator and the commutator $[T_1, U_1]$ is compact.
\item $\norm{U_2} < \epsilon. $

\end{enumerate}

Then $T+U$ is a left  pseudo semi-B-Fredholm operator and $ \epsilon$ can be chosen so that $ \alpha(T_2+U_2) \leq \alpha(T_2)$ and if  $ \beta(T_2) = \infty , $  $ \beta(T_2+U_2)= \infty.$

\eet

\bp Let $T$ be a   left  pseudo-semi-B-Fredholm operator and let $P$  be an idempotent associated to $T$  in $\mathcal{C}(X).$ Then $ T_2$ is a left semi-Fredholm operator. From \cite[Theorem 4.2.1]{CPY}, there exists $\epsilon_1 > 0,$ such that if $ V\in L(X_2) $ satisfies  $\norm{V} < \epsilon_1, $ then
$ T_2+ V$ is an upper semi-Fredholm operator , the nullity of $T_2 + V$  and the deficiency of $T_2 +V$  satisfies:

%\begin{itemize}

 $\bullet \alpha(T_2+V)  \leq \alpha(T_2),  $

$\bullet$ If  $ \beta(T_2) = \infty , $ then  $ \beta(T_2+V)= \infty.$
%\end{itemize}

Let $ U \in L(X) $ such that    $\Pi(P)$ commutes with $ \Pi(U)$  and $\norm{U_2} < \epsilon_1. $
Then
$ T_2+ U_2$ is an upper semi-Fredholm operator , the nullity of $T_2 + U_2$  and the deficiency of $T_2 +U_2$  satisfies:

%\begin{itemize}

 $\bullet$ $\alpha(T_2+ U_2)  \leq \alpha(T_2),  $ (*)

$\bullet$ If  $ \beta(T_2) = \infty , $ then  $ \beta(T_2+ U_2)= \infty.$  (**)
%\end{itemize}

Again as  $ T_2$ is a left semi-Fredholm operator, its image   in the Calkin algebra $\mathcal{C}(X_2)$ is left invertible.   As the set of left invertible elements in a Banach algebra is open, there exists  $\epsilon_2 > 0,$ such that if $ U \in L(X) $ satisfies $\Pi(P)$ commutes with $ \Pi(U)$  and  $\norm{U_2} < \epsilon_2,$  then
$(T+ U)_2=  T_2+ U_2 $ is a left  semi-Fredholm operator.

Let $ \epsilon< min( \epsilon_1, \epsilon_2),$  then  if \, $ U \in L(X) $ satisfies     $\Pi(P)$ commutes with $ \Pi(U),$   $\norm{U_2} < \epsilon,$
 then $ (T + U)_2= T_2 + U_2 $ is a left semi-Fredholm operator and the conditions (*) and (**) are satisfied.\\
 As $\Pi(P)$ commutes with $ \Pi(U),$ then $\Pi(P)$ commutes with $ \Pi(T+U).$ Moreover as $T_1$ and $U_1$ are both  Riesz operators and the commutator $[T_1, U_1]$ is compact, then $T_1+ U_1 $ is a Riesz operator. So $ T+U= (T_1+U_1) \oplus ( T_2 + U_2) +K= (T+U)_1 \oplus ( T + U)_2 +K,$ where $K \in L(X)$ is a compact operator.  From \cite[Theorem 3.5]{HB}, it follows that $T+U$ is a left pseudo semi-B-Fredholm operator.  \ep

In a similar way to Theorem \ref{left-perturbation},  using \cite[Theorem 4.2.2]{CPY}, the fact that the set of right invertible elements in a Banach algebra is open and  \cite[Theorem 3.5]{HB}, we can prove by the same method the following result.

\bet \label{right-perturbation} Let $T$ be a   right  pseudo semi-B-Fredholm operator and let $P \in L(X)$  be an idempotent associated to $T$ in $\mathcal{C}(X).$  Then there exists an $ \epsilon > 0,$ such that if $U \in L(X)$  is such that:

\begin{enumerate}
\item  $\Pi(P)$ commutes with $ \Pi(U),$
\item  $U_1$ is a Riesz operator commuting with $T_1,$
\item $\norm{U_2} < \epsilon. $
\end{enumerate}

Then $T+U$ is a  right pseudo semi-B-Fredholm operator and $ \epsilon$ can be chosen so that $ \beta(T_2+U_2) \leq \beta(T_2)$ and if  $ \alpha(T_2) = \infty , $  $ \alpha(T_2+U_2)= \infty.$

\eet

Theorem \ref{left-perturbation} and Theorem \ref{right-perturbation} give us immediately the following result.

\bet \label{Fred-perturbation} Let $T$ be a  pseudo B-Fredholm operator and let  $P \in L(X)$  be an idempotent associated to $T$ in $\mathcal{C}(X).$   Then there exists an $ \epsilon > 0$ such that if \, $U \in L(X)$  satisfies:

\begin{enumerate}

\item  $\Pi(P)$ commutes with $ \Pi(U),$
\item  $U_1$ is a Riesz operator and the commutator $[ T_1, U_1]$ is compact.
\item $\norm{U_2} < \epsilon. $

\end{enumerate}

Then $T+U$ is a  pseudo B-Fredholm operator and   \textit{\textbf{{ind}}}($T+U$)= \textit{\textbf{{ind}}}($T$).

\eet

\bp  From Theorem \ref{left-perturbation} and Theorem \ref{right-perturbation}, there exists  $ \epsilon_1 > 0,$ such that if $ U \in L(X)$  satisfies
$\Pi(P)$ commutes with $ \Pi(U),$   $\norm{U_2} < \epsilon, $   $U_1$ is a Riesz operator and the commutator $[ T_1, U_1]$ is compact,
then  $T+U$ is a  pseudo B-Fredholm operator. As $ T_2$ is a Fredholm operator, there exists $  \epsilon_2 > 0,$ such that if  $ V\in L(X_2) $ and  $\norm{V}< \epsilon_2,$
then  $ T_2 + V$ is a Fredholm operator and  $ ind(T_2 + V) = ind(T_2)$

Let $ \epsilon< min( \epsilon_1, \epsilon_2).$ So if   $\norm{U_2} < \epsilon, $ then $(T+U)_2= T_2+ U_2$ is a Fredholm operator and $ind( (T+U)_2 )= ind( T_2).$  Thus \textit{\textbf{{ind}}}($T+U$)= \textit{\textbf{{ind}}}($T$). \ep

\bremark  When $P=0,$  $\Pi(P)=0$ commutes with all elements of $\mathcal{C}(X)$ and  we retrieve in the previous results, the usual properties of the Fredholm index.
\eremark

 \noindent  Alaa Hamdan\\
\noindent   Dubai \\
\noindent   United Arab Emirates \\
\noindent     aa.hamdan@outlook.com \\

\noindent  Mohammed Berkani\\
\noindent Honorary member of LIAB \\
\noindent Department of Mathematics\\
\noindent Science Faculty of Oujda \\
\noindent Mohammed I University\\
\noindent Morocco\\
\noindent berkanimoha@gmail.com


\begin{thebibliography}{5}




%\bibitem{AA} Y.A. Abramovich, C.D.Aliprantis, \textit {Problems in operator theory,} Graduate Texts in Mathematic, Volume 51,
%American Mathematical Society, Providence, Rhode Island.

\bibitem{AI}P. Aiena, \textit{ Fredholm and Local Spectral Theory with Applications Multipliers},
Kluwer Academic Publishers, 2004. London.

\bibitem{AZOUZ} Z. Aznay, A. Ouahab, H. Zariouh,
\textit{ On the index of pseudo B-Fredholm operator,}
Monatshefte f$\ddot{u}$r Mathematik volume 200, pages 195-207 (2023).


\bibitem{BAR} B. Barnes,  \textit{ Essential Spectra in a Banach Algebra Applied to Linear
  Operators,} Proceedings of the Royal Irish Academy. Section A: Mathematical
and Physical Sciences Vol. 90A, No. 1 (1990), pp. 73-82.


%\bibitem{BMW}  B. Barnes, G. J. Murphy, M.R.F. Smyth, T.T. West,
% \textit{Riesz and Fredholm theory in Banach algebras}, Pitman Publishing
%Inc, 1982. London.

\bibitem{P7}M. Berkani,
 \textit{ On a class of quasi-Fredholm operators,} Integr. Equ. Oper.
Theory, 34 (1999), 244-249.

\bibitem{P10} M. Berkani and M. Sarih,  \textit{  An Atkinson
type theorem for $B$-Fredholm operators,}  Studia Math. {148}
(2001), 251--257.

\bibitem{P45} M. Berkani, S. \v{C}. \v{Z}ivkovi\'{c}-Zlatanovi\'{c}, \textit{ Pseudo-B-Fredholm operators, poles of
the resolvent and mean convergence in the Calkin algebra,}
 Filomat 33 (2019), no. 11, 3351-3359.


\bibitem{P48}M. Berkani, \textit{ A new approach in index theory,}  Mathematische Zeitschrift,  volume 298,no. 3-4, 943-951 (2021).

\bibitem{P49}{M. Berkani} \textit{ Continuous Fredholm Theory, Regularities and Semiregularities.} Complex Anal. Oper. Theory 15, 105 (2021).
https://doi.org/10.1007/s11785-021-01151-1.


\bibitem {P50} M. Berkani,   \textit{On the  p-invertibility  in   rings and Banach algebras},
Aequat. Math. 96, 1125-1140 (2022). https://doi.org/10.1007/s00010-022-00892-4

%\bibitem{BS} P. S. Bourdon, J. H. Shapiro,     \textit{Riesz Composition Operators},
%Pacific Journal of Mathematics
%Vol. 181,  No. 2 1997

\bibitem{BO}E. Boasso,\textit{Isolated spectral points and Koliha-Drazin invertible
elements in quotient Banach algebras and homomorphisms ranges,}
Mathematical Proceedings of the Royal Irish Academy 115 (2), 2015,
1-15.


\bibitem{CPY}  S.R. Caradus, W.E. Pfaffenberger, Bertram Yood,
\textit{Calkin Algebras and Algebras of Operators on Banach Spaces},
Lecture Notes in Pure and Applied Mathematics,
MARCEL DEKKER, INC. New York 1974.

\bibitem {DR} { M. P. Drazin,} \textit{ Pseudoinverse in
associative rings and semigroups,}  Amer. Math. Monthly 65 (1958),
506-514.

\bibitem {HA} R. E. Harte, \textit{  On quasinilpotents in rings.} PanAm.
Math. J. 1 (1991), 10-16.

\bibitem{Koliha} J. J. Koliha, \textit{A generalized Drazin inverse}, Glasg. Math. J. 38 (1996),
367-381.
.
\bibitem{KP} J. J. Koliha, P. Patricio, \textit{Elements of rings
 with equal spectral idempotents},
J. Austral. Math. Soc. 72 (2002), 137-152.

\bibitem{HB} A. Hamdan, M. Berkani, \textit{On Classes of Fredholm Type Operators},
https://arxiv.org/abs/2401.01061.



%\bibitem{TAKAM}  A. Tajmouati, M. Karmouni, Mbark Abkari,
%\textit{Pseudo semi B-Fredholm and Generalized Drazin invertible operators Through Localized SVEP,}
%Italian Journal of Pure and Applied Mathematics  37-2017 (301-314).




%\bibitem{P11} M. Berkani and M. Sarih, {\it On semi B-Fredholm
%operators}, Glasgow Math. J. 43 (2001), p. 457-465.

%\bibitem{P42}{M. Berkani},
%{\it A trace formula for the index of B-Fredholm operators,}
%http://arxiv.org/abs/1609.01695


%\bibitem{P13} {M. Berkani}, {\it B-Weyl spectrum and poles of the
%resolvent}, J. Math. Anal. App. 272 (2002), 596-603.

%\bibitem{BEL}O. Bel Hadj Fredj, On the poles of the
%resolvent in Calkin algebras, Proc. Amer. Math. Soc. 135 (2007),
%2229-2234.

%\bibitem{P8}{M. Berkani},
%{\em Restriction of an operator to the range of its powers\/},
%Studia Mathematica, 140 (2)\, (2000), 163-174.

%\bibitem{P12}M. Berkani,
%Index of B-Fredholm operators and generalization of a Weyl
%Theorem,  Proc. Amer. Math. Soc., 130 (6),\, (2002), 1717-1723.

%\bibitem{P33}{M. Berkani},
%{\em On the B-Fredholm Alternative\/}, Mediterr. J. Math. ,
%10(3),\, {2013},\,1487-1496.

%\bibitem{P10} M. Berkani and M. Sarih,    An Atkinson
%type theorem for $B$-Fredholm operators,  Studia Math. {148}
%(2001), 251--257.

%\bibitem {BMM} M. Benharrat, K.   Miloud Hocine, B.
%Messirdi, { \em Left and right generalized Drazin invertible
%operators and local spectral theory.} Proyecciones 38 (2019), no.
%5, 897-919.

%\bibitem{BO}E. Boasso, Isolated spectral points and Koliha-Drazin invertible
%elements in quotient Banach algebras and homomorphisms ranges.
%Mathematical Proceedings of the Royal Irish Academy 115 (2), 2015,
%1-15.

%\bibitem{BOJ} {M.T. Boedihardjo,} {\it Topics in functional analysis,}
%Dissertation, Texas AM University, 2016.

%\bibitem{BOWJ} {M.T. Boedihardjo; W.B. Johnson,} {On mean ergodic
%convergence in the Calkin algebras.} Proc. Amer. Math. Soc. 143
%(2015), no. 6, 2451�2457.

%\bibitem{BD} { F.F. Bonsall, J.Duncan,}  {\it Complete Normed Algebras,}
%Springer-Verlag, 1973. Berlin.

%\bibitem{CA}{S.R. Caradus,}  { \it Operator theory of the
%pseudo-inverse,} Queen's papers in Pure and Appl. Math. No.
%38(1974).

%\bibitem{BUR} L. Burlando,  A generalization of the uniform ergodic theorem to
%poles of arbitrary order,   Studia Mathematica 122 -1, 1997,75-98.

%\bibitem{CZ} M.D. Cvetkovi\'c and S. \v{C}. \v{Z}ivkovi\'c-Zlatanovi\'c,
%Generalized Kato decomposition and essential  spectra, Compl.
%Anal. Oper. Th., (2017),  11(6),  1425-1449.


%\bibitem{CBS} { M.D. Cvetkovic, E. Boasso, S.C.
%Zivkovic-Zlatanovic,} { \it Generalized B-Fredholm Banach algebra
%elements,} Mediterr. J. Math. (2016) 13: 3729.
%doi:10.1007/s00009-016-0711-y



%\bibitem{DU} { Dunford, N. \& Schwartz, J.} { Linear operators,
%Part 1;}  Wiley Inter-science, New York (1971).\\


%\bibitem {FLT} D.E. Ferreyra, M.Lattanzi, F.E. Levis, N. Thome,
% { \em Left and right generalized Drazin invertible operators on Banach spaces and
% applications.}
 %Oper. Matrices 13 (2019), no. 3, 569-583.

%\bibitem{GR}  S., Grabiner Uniform ascent and
%descent of bounded operators, J. Math. Soc. Japan 34 , no 2 ,
%(1982), 317-337.


%\bibitem{GRZ} S. Grabiner, J. Zemanek, Ascent, descent, and ergodic properties of linear operators
%  J. Operator Theory 48(2002), 69-81

%\bibitem{TAJ}A. Tajmouati, M. Karmouni, On pseudo B-Weyl and pseudo
%B-Fredholm operators, Int. J. Pure Appl. Math. 108, 513-522.
%(2016).

%\bibitem{ZZ} H. Zariouh, H. Zguitti, On pseudo B-Weyl operators
%and generalized Drazin invertibility for operator matrices, Linear
%and Multilinear Algebra 64 (7), 1245-1257, 2016.


\end{thebibliography}
\end{document}